\documentclass[12pt,reqno]{amsart}

\usepackage{fullpage,amsmath,amssymb,amsfonts,mathrsfs,bm,graphicx,hyperref}

\def\bbone{{\mathchoice {\rm 1\mskip-4mu l} {\rm 1\mskip-4mu l}
{\rm 1\mskip-4.5mu l} {\rm 1\mskip-5mu l}}}

\newcommand{\ux}{\mathbf{x}}

\newtheorem{theorem}{Theorem}[section]

\newtheorem{proposition}{Proposition}[section]

\begin{document}

\author{Abdelmalek Abdesselam}
\address{Abdelmalek Abdesselam, Department of Mathematics,
P. O. Box 400137,
University of Virginia,
Charlottesville, VA 22904-4137, USA}
\email{malek@virginia.edu}

\title{A local injective proof of log-concavity for increasing spanning forests}

\begin{abstract}
We give an explicit combinatorial proof of a weighted version of strong log-concavity for the generating polynomial of increasing spanning forests of a finite simple graph equipped with a total ordering of the vertices. In contrast to similar proofs in the literature, our injection is local in the sense that it proceeds by moving a single edge from one forest to the other. In the particular case of the complete graph, this gives a new combinatorial proof of log-concavity of unsigned Stirling numbers of the first kind where a pair of permutations is transformed into a new pair by breaking a single cycle in the first permutation and gluing two cycles in the second permutation, while all the other spectator cycles are left untouched.
\end{abstract}

\maketitle

\section{Introduction}

Let $G$ a finite simple graph with vertex set $V=\{v_1,v_2,\ldots,v_n\}$
and edge set $E$. We assume we are given a total ordering of the vertices $v_1<v_2<\cdots<v_n$.
A spanning forest $\mathscr{F}$ is a subgraph of $G$ with the same vertex set and an edge set given by a subset of $E$ which contains no circuit. Recall that a circuit is a set of edges (seen as unordered pairs of vertices) of the form $\{\{v_{i_1},v_{i_2}\},\{v_{i_2},v_{i_3}\},\ldots,
\{v_{i_{p-1}},v_{i_p}\},\{v_{i_p},v_{i_1}\}\}$, with $p\ge 3$ and where all the vertices $v_{i_1},\ldots,v_{i_p}$ are distinct.
Such a forest $\mathscr{F}$ is thus made of a disjoint collection of trees which will, by convention, be rooted at their minimal vertex. Such a forest is called {\em increasing} if for every such tree, the vertices (or labels) increase along each path from the root to a leaf.
For simplicity, we will restrict to $V=[n]:=\{1,2,\ldots,n\}$ with the usual ordering. We will also identify a forest $\mathscr{F}$ with the corresponding subset of edges inside $E$. We will also switch gears and view an edge $e$ as an ordered pair of vertices $(i,j)$ with $i<j$.
For example, with $G$ being the complete graph on nine vertices,
\begin{equation}
\mathscr{F}_1=\{
(1,2),(1,4),(4,7),(4,9),(3,5),(3,6),(6,8)
\}
\label{forest1}
\end{equation}
is an increasing forest, whereas
\[
\mathscr{F}_2=\{
(1,2),(1,8),(7,8),(8,9),(3,5),(3,6),(4,6)
\}
\]
is not, because of the drops occurring in the paths from $1$ to $7$ and from $3$ to $4$.

\smallskip
Increasing forests and associated generating polynomials have been studied recently in~\cite{HallamS,HallamMS}. Increasing forests have also been investigated in~\cite{AbreuN} in connection to LLT polynomials. Given a collection of indeterminates $\ux=(x_e)_{e\in E}$
and for $0\le k\le n$, one can define the polynomial
\[
a_k(\ux)=\sum_{\mathscr{F}\in {\rm IF}_k}\ \prod_{e\in\mathscr{F}} x_e
\]
where the sum is over the set ${\rm IF}_k$ of all increasing forests with exactly $k$ connected components.
Given an additional variable $t$, such polynomials can be packaged into the generating polynomial
\[
{\rm ISF}(\ux,t)=\sum_{k=0}^{n}a_k(\ux)t^k\ .
\]
Of course, when $n\ge 1$ the first term $a_0(\ux)$ vanishes and all these polynomials depend on the graph $G$ and the chosen total ordering. 
Our main concern in this article will be the log-concavity of the sequence $a_k(\ux)$, $0\le k\le n$.

We now recall some basic terminology about log-concavity. The reader can glean much more information from the reviews~\cite{Stanley,Brenti,Branden}.
A finite sequence $(a_0,a_1,\ldots,a_n)$ of real numbers is called log-concave if
$a_{p-1}a_{p+1}\le a_p^2$ whenever $0<p<n$. 
It is called strongly log-concave if $a_{p-1}a_{q+1}\le a_{p}a_q$ whenever $0<p\le q<n$. Such a sequence is said to be without internal zeros if there are no indices $0\le p<q<r\le n$ with $a_p\neq 0$, $a_q=0$ and $a_r\neq 0$.
Of course, log-concavity follows from strong log-concavity as the $p=q$ particular case.
However, when the sequence is made of nonnegative numbers and has no internal zeros, the two properties are easily seen to be equivalent.
Following~\cite{Sagan2,Sagan3}, we will define a weighted version of these two notions for sequences of polynomials like $(a_0(\ux),a_1(\ux),\ldots,a_n(\ux))$ in a polynomial algebra $\mathbb{R}[\ux]$. We will say that the sequence is $\ux$-log-concave if, whenever
$0<p<n$, it holds that the polynomial $a_p(\ux)^2-a_{p-1}(\ux)a_{p+1}(\ux)$ only has nonnegative coefficients. Likewise, we will say that the sequence is strongly $\ux$-log-concave if, whenever
$0<p\le q<n$, it holds that the polynomial $a_p(\ux)a_q(\ux)-a_{p-1}(\ux)a_{q+1}(\ux)$ only has nonnegative coefficients. In general, for the weighted case, strong $\ux$-log-concavity is a strictly stronger property than $\ux$-log-concavity. Although we did not find it explicitly stated in the literature, one has the following $\ux$-log-concavity result.

\begin{theorem}
\label{easythm}
The sequence $(a_0(\ux),a_1(\ux),\ldots,a_n(\ux))$ defined above by the enumeration of increasing spanning forests is strongly $\ux$-log-concave. 
\end{theorem}

Indeed, that Hallam et al.~\cite[Theorem 2.2]{HallamMS} have shown the factorization
\[
{\rm ISF}(\ux,t)=\prod_{j=1}^{n}\left(t+
\sum_{i=1}^{j-1}\bbone\{(i,j)\in E\}\ x_{(i,j)}\right)
\]
where $\bbone\{\cdots\}$ denotes the indicator function of the logical condition between braces.
Therefore, Theorem \ref{easythm} follows immediately from this factorization together with~\cite[Proposition 2.3]{Sagan3}. 
Also note that, if one specializes the weights to real values, then the log-concavity of the $a_k$ sequence is immediate since they form the coefficients of a polynomial with only real roots. One even gets ultra log-concavity, namely, the log-concavity of the sequence $a_k/\binom{n}{k}$, see~\cite[Theorem 2]{Stanley}.
In this article we will provide a new
proof of Theorem \ref{easythm} which is based on an explicit weight-preserving combinatorial injection, and which is thematically similar to many log-concavity proofs in the literature, e.g.,~\cite{Gasharov,Hamidoune,Krattenthaler,Leroux,Sagan1,Sagan2,Sagan3}. The
key difference however, is what we, for lack of a better word, call {\em locality}.

Let $0\le k<\ell\le n$, and let $\Psi:{\rm IF}_k\times {\rm IF}_{\ell}
\rightarrow {\rm IF}_{k+1}\times{\rm IF}_{l-1}$
be an injective map. We will denote the image by $\Psi$ of a pair of increasing forests $(\mathscr{A},\mathscr{B})$ simply by $(\mathscr{A}',\mathscr{B}')$.
We will say that $\Psi$ is {\em local} if for all $(\mathscr{A},\mathscr{B})\in {\rm IF}_k\times {\rm IF}_{\ell}$, there exists an edge $e\in\mathscr{A}\backslash\mathscr{B}$, such that
\[
\mathscr{A}'=\mathscr{A}\backslash\{e\}\qquad \text{and}\qquad
\mathscr{B}'=\mathscr{B}\cup\{e\}\ .
\]
We will refer to such an edge as a movable edge.
Note that locality forces weight preservation, i.e.,
\[
\left(\prod_{e\in\mathscr{A}'}x_e\right)
\times\left(\prod_{e\in\mathscr{B}'}x_e\right)
=
\left(\prod_{e\in\mathscr{A}}x_e\right)
\times\left(\prod_{e\in\mathscr{B}}x_e\right)\ .
\]
Now our main result is as follows.

\begin{theorem}
\label{mainthm}
For all $0\le k<\ell\le n$, there exists an explicit local injective map
$\Psi:{\rm IF}_k\times {\rm IF}_{\ell}
\rightarrow {\rm IF}_{k+1}\times{\rm IF}_{l-1}$\ .
\end{theorem}

Making the change of indices $p=k+1$ and $q=\ell-1$, and thanks to weight preservation, one easily sees that Theorem \ref{easythm} follows from Theorem \ref{mainthm}. Note that we now include the bijective case $q=p-1$ or $\ell=k+1$ which is silly from the point of view of strong log-concavity.
A trivial bijection ${\rm IF}_k\times {\rm IF}_{k+1}
\rightarrow {\rm IF}_{k+1}\times{\rm IF}_{k}$ that does the job is the forest swap given by
$(\mathscr{A}',\mathscr{B}'):=(\mathscr{B},\mathscr{A})$.
Locality makes the existence of such a bijection less trivial since one is only allowed to transfer a single edge.
In the particular case of the complete graph and the Stirling numbers of the first kind, 
to be discussed in \S\ref{Stirlingsec}, our injection for pairs of permutations $(\sigma,\tau)\mapsto(\sigma',\tau')$ is such that $\sigma'$ is obtained from $\sigma$ by breaking a single cycle, whereas $\tau'$ is obtained by gluing two cycles of $\tau$. We call this a {\em local move} as opposed to a global transformation where the resulting premutations can look very different from the original ones. This global flavor was a feature of, e.g.,
previous injections by Sagan~\cite{Sagan1,Sagan3} which use the powerful Lindstr\"om-Gessel-Viennot machinery of lattice paths~\cite{Lindstrom,GesselV}.

\section{The explicit injection}

\subsection{The construction}
An important tool for the construction of our local injection is provided by the next proposition. We will use absolute values to denote the cardinality of finite sets, and for $0\le k\le n$, we will denote by ${\rm S}_{n,k}$ the set of subsets of $[n]$ which have $k$ elements.

\begin{proposition}\label{toolprop}
For all $n\ge 0$ and $0\le k<\frac{n}{2}$, there exists an injection
$\Phi_{n,k}:{\rm S}_{n,k}\rightarrow {\rm S}_{n,k+1}$ such that for all $k$-subset $X$, we have the inclusion $X\subset\Phi_{n,k}(X)$.
\end{proposition}

There are many ways to produce such injections, for instance, the bracketing algorithm in~\cite{GreeneK} which is an elegant reformulation of earlier constructions~\cite{deBruijnEtal,Aigner} (see~\cite{Pak,Zeilberger} for insightful reviews on this topic). Let us also sketch another approach which echoes the remarks at the end of the previous section about emphasizing the ``apparently silly bijective case''. First construct all the bijections
$\Phi_{2k+1,k}$, for all $k$. This amounts to chosing a perfect matching in the bipartite Kneser graph $H(2k+1,k)$, i.e., the middle two levels of the Boolean lattice of $[2k+1]$. There are many matchings to choose from (see~\cite{KiersteadT} and the more recent~\cite{Jin}) and their study was motivated by the middle levels conjecture (see~\cite{GregorMN} and references therein).
Then $\Phi_{n,k}$ can be constructed recursively, or by induction on $n$ as follows. If $k<\frac{n-1}{2}$, then let
\[
\Phi_{n,k}(X):=
\left\{
\begin{array}{ll}
\Phi_{n-1,k-1}(X) & \text{if}\ n\notin X\ , \\
\{n\}\cup\Phi_{n-1,k-1}(X\backslash\{n\}) & \text{if}\ n\in X\ .
\end{array}
\right.
\]
The remaining case $\frac{n-1}{2}\le k<\frac{n}{2}$, which can only happen if $n=2k+1$, has been taken care of beforehand.
More generally for a finite set $Y$ instead of $[n]$, and for $0\le k< |Y|/2$
we will denote by $\Phi_{Y,k}$ a choice of injection from $k$-subsets to $(k+1)$-subsets of $Y$ satisfying the same constraint as in Proposition \ref{toolprop}.

We now proceed with the definition of our injection $\Psi$. 
For each subset $Y\subset [n]$ and each $k$, $0\le k<|Y|/2$, we prepare the maps $\Phi_{Y,k}$ discussed above.
For $\mathscr{A}$ an increasing forest, we will use the notation $m(\mathscr{A})$
for the set of minima of connected components of $\mathscr{A}$. For $(\mathscr{A},\mathscr{B})\in {\rm IF}_k\times {\rm IF}_{\ell}$ we therefore have $|m(\mathscr{A})|=k$
and $|m(\mathscr{B})|=\ell$. From the hypothesis $k<\ell$, we see that the symmetric difference
\[
m(\mathscr{A})\Delta m(\mathscr{B})=\left(m(\mathscr{A})\backslash m(\mathscr{B})\right)
\cup \left(m(\mathscr{B})\backslash m(\mathscr{A})\right)
\]
satisfies
\[
|m(\mathscr{A})\backslash m(\mathscr{B})|<\frac{1}{2}\ |m(\mathscr{A})\Delta m(\mathscr{B})|\ .
\]
As a result, we have a uniquely defined element $j\in [n]$ for which
\[
\{j\}=\Phi_{m(\mathscr{A})\Delta m(\mathscr{B}),|m(\mathscr{A})\backslash m(\mathscr{B})|}
\left(m(\mathscr{A})\backslash m(\mathscr{B})\right)\backslash \left(m(\mathscr{A})\backslash m(\mathscr{B})\right)\ .
\]
Let $A\subset [n]$ be the set of vertices corresponding to the connected component of the forest $\mathscr{A}$ which contains $j$. Likewise, let $B\subset [n]$ be the set of vertices corresponding to the connected component of the forest $\mathscr{B}$ which contains $j$. By definition, $j$ is the minimal element of $B$ but it is different from the minimal element $i_0$ of $A$. Consider the path in the forest $\mathscr{A}$ going from $i_0$ to $j$ and let $e=(i,j)$ be the last edge on that path. We claim that $e$ is a movable edge, i.e., that $e\in\mathscr{A}\backslash\mathscr{B}$ and, letting
\begin{equation}
\mathscr{A}':=\mathscr{A}\backslash\{e\}\qquad \text{and}\qquad
\mathscr{B}':=\mathscr{B}\cup\{e\}
\label{primedef}
\end{equation}
we have that both $\mathscr{A}'$, $\mathscr{B}'$ are increasing forests.
Clearly, removing an edge from an increasing forest produces an increasing forest so we only need explain the property for $\mathscr{B}'$. Since $\mathscr{A}$ is increasing $i_0\le i<j=\min(B)$ so $i$ belongs to a different connected component $\widetilde{B}$ of $\mathscr{B}$. Since $e$ joins different connected components, adding it to $\mathscr{B}$ still gives a forest. Let $\widetilde{j}$ be the minimal element of $\widetilde{B}$. It will also be the minimal element of the newly formed connected component $B\cup\widetilde{B}$ of the forest $\mathscr{B}'$.
It is easy to see that paths from the root $\widetilde{j}$ to a leaf are increasing. If the leaf is in $\widetilde{B}$, it is so because it already was an increasing path in $\mathscr{B}$. If the leaf $b$ is in $B$, then the path must be of the form $\widetilde{j},\ldots,i,j,\ldots,b$ where the portion $\widetilde{j},\ldots,i$ is in $\widetilde{B}$ and the portion $j,\ldots,b$ is in $B$. Both portions increase because $\mathscr{B}$ is an increasing forest, while at the transition $i<j$ by construction.
Of course the new increasing forests $\mathscr{A}',\mathscr{B}'$ respectively have $k+1$ and $\ell-1$ connected components. Hence,
we now have a well defined map $\Psi$ from ${\rm IF}_k\times {\rm IF}_{\ell}$ to
${\rm IF}_{k+1}\times {\rm IF}_{\ell-1}$.

\subsection{Proof of the injective property}

Suppose two pairs of forests $(\mathscr{A}_1,\mathscr{B}_1)$
and $(\mathscr{A}_2,\mathscr{B}_2)$ are sent to the same $(\mathscr{A}',\mathscr{B}')$.
Let $j_1$ be the component minimum, and $e_1$ the movable edge, as in the previous section, which appeared in the construction of  $(\mathscr{A}',\mathscr{B}')$ from 
$(\mathscr{A}_1,\mathscr{B}_1)$. Likewise, let $j_2$ be the component minimum, and $e_2$ the movable edge which 
appeared in the construction
of  $(\mathscr{A}',\mathscr{B}')$ from $(\mathscr{A}_2,\mathscr{B}_2)$.
It is easy to track the changes in the sets of component minima through the map $\Psi$ which results in the relations
\[
m(\mathscr{A}')=m(\mathscr{A}_1)\cup\{j_1\}=m(\mathscr{A}_2)\cup\{j_2\}
\]
and
\[
m(\mathscr{B}')=m(\mathscr{B}_1)\backslash\{j_1\}=m(\mathscr{B}_2)\backslash\{j_2\}
\]
with the knowledge that $j_1\in m(\mathscr{B}_1)\backslash m(\mathscr{A}_1)$ and $j_2\in m(\mathscr{B}_2)\backslash m(\mathscr{A}_2)$.
We thus have preservation of unions:
\[
m(\mathscr{A}')\cup m(\mathscr{B}')=
m(\mathscr{A}_1)\cup m(\mathscr{B}_1)=
m(\mathscr{A}_2)\cup m(\mathscr{B}_2)
\]
and similarly for intersections and symmetric differences (i.e., the same equalities as above hold
with $\cap$, $\Delta$, instead of $\cup$).
From $|m(\mathscr{A}_1)|=|m(
\mathscr{A}_2)|=k$, and denoting by $r$ the cardinality of $m(\mathscr{A}_1)\cap m(\mathscr{B}_1)=
m(\mathscr{A}_2)\cap m(\mathscr{B}_2)$, we get
\[
|m(\mathscr{A}_1)\backslash m(\mathscr{B}_1)|
=|m(\mathscr{A}_2)\backslash m(\mathscr{B}_2)|=k-r\ .
\]
Therefore
\begin{eqnarray*}
m(\mathscr{A}')\backslash m(\mathscr{B}') & = & \left(
m(\mathscr{A}_1)\backslash m(\mathscr{B}_1)
\right)\cup\{j_1\} \\
 & = & \Phi_{m(\mathscr{A}_1)\Delta m(\mathscr{B}_1),
 |m(\mathscr{A}_1)\backslash m(\mathscr{B}_1)|}
\left(m(\mathscr{A}_1)\backslash m(\mathscr{B}_1)\right) \\
 & = & \Phi_{m(\mathscr{A}')\Delta m(\mathscr{B}'),k-r}
\left(m(\mathscr{A}_1)\backslash m(\mathscr{B}_1)\right)
\end{eqnarray*}
and similarly for the other pair. Therefore
\[
\Phi_{m(\mathscr{A}')\Delta m(\mathscr{B}'),k-r}
\left(m(\mathscr{A}_1)\backslash m(\mathscr{B}_1)\right)=
\Phi_{m(\mathscr{A}')\Delta m(\mathscr{B}'),k-r}
\left(m(\mathscr{A}_2)\backslash m(\mathscr{B}_2)\right)
\]
and the injectivity of the $\Phi$ maps forces $m(\mathscr{A}_1)\backslash m(\mathscr{B}_1)=
m(\mathscr{A}_2)\backslash m(\mathscr{B}_2)$. Since the unions and intersections also agree, we immediately obtain $m(\mathscr{A}_1)=m(\mathscr{A}_2)$ and $m(\mathscr{B}_1)=m(\mathscr{B}_2)$.
As a result
\[
\{j_1\}=m(\mathscr{A}')\backslash m(\mathscr{A}_1)
=m(\mathscr{A}')\backslash m(\mathscr{A}_2)=\{j_2\}\ .
\]
The equality $j_1=j_2$, forces that of the unique edges in $\mathscr{B}'$ of the form $(i_1,j_1)$ and $(i_2,j_2)$ with $i_1<j_1$ and $i_2<j_2$, namely the equality $e_1=e_2$. 
Finally, the definition (\ref{primedef})
gives the desired conclusion $(\mathscr{A}_1,\mathscr{B}_1)=(\mathscr{A}_2,\mathscr{B}_2)$.

\section{Concluding remarks}

\subsection{Log-concavity for elementary symmetric functions}\label{subsetsec}
Note that Proposition \ref{toolprop} gives a combinatorial proof of unimodality of binomial coefficients. One can adapt the construction of our $\Psi$ to this simpler setting and obtain a new derivation of the log-concavity of binomial coefficients. Indeed, one can construct a map $S_{n,k}\times S_{n,\ell}\rightarrow S_{n,k+1}\times S_{n,\ell-1}$, $(X,Y)\mapsto (X',Y')$
by letting
\begin{eqnarray*}
\{i\} & := & \Phi_{X\Delta Y,|X\backslash Y|}(X\backslash Y)\backslash (X\backslash Y) \\
X' & := & X\cup\{i\} \\
Y' & := & Y\backslash\{i\}
\end{eqnarray*}
just as in (\ref{primedef}).
This map is injective and also weight-preserving. Now  taking $\ux=(x_1,\ldots,x_n)$ and weighting subsets by the product of corresponding $x$ variables, we obtain a combinatorial proof of the strong $\ux$-log-concavity of the elementary symmetric functions $e_k(\ux)$ as in~\cite{Sagan3}. Of course, by the Jacobi-Trudy identity, we know more, namely, that the relevant polynomial difference is not only positive coefficient-wise in the $\ux$, but is also Schur-positive.
This is one of the simplest instances of a very general phenomenon as investigated, e.g., in~\cite{LamPP}.

\subsection{Stirling numbers of the first kind}\label{Stirlingsec}

We will denote by $s(n,k)$ the Stirling numbers of the first kind and by $c(n,k):=|s(n,k)|$ their unsigned version. As is well known, $c(n,k)$ counts the number of permutations of $n$ elements with exactly $k$ cycles. It is easy to construct a bijection between such permutations and increasing forests in ${\rm IF}_k$, where $G$ is the complete graph on $[n]$. Let $\mathscr{F}\in{\rm IF}_k$, and turn it into a collection of rooted planar trees as follows. For each tree, draw the root (again the minimal vertex) on the left and the leafs on the right. For a vertex $v$ with children $w_1<\cdots<w_p$ draw them, in this order, from top to bottom. Then consider a path around the tree starting and ending at the root, going in the counterclockwise direction. Intuitively, one could think of drawing the contour of one's hand placed flat on a table. Next, make the cycle $(i_1 i_2\cdots i_q)$ where $i_1$ is the root and the other vertices are listed in the order in which they appear as one turns around the tree. Finally, let $\sigma$ be the permutation obtained as the product of such cycles, over all trees in the forest $\mathscr{F}$. For example, the permutation corresponding to the forest $\mathscr{F}_1$ in (\ref{forest1}), in cycle notation, is
\[
\sigma=(14972)(3685)\ .
\]
Using this bijection and our map $\Psi$, we thus obtain an injective proof of log-concavity of the sequence $(c(n,0),\ldots,c(n,n))$ which is different than the constructions in~\cite{Sagan1,Sagan3}. In particular our map sends a pair of permutations $(\sigma,\tau)$ to a new one $(\sigma',\tau')$ where $\sigma'$ is obtained by breaking a cycle of $\sigma$ in two, and $\tau'$ is obtained by gluing two cycles of $\tau$, while all other cycles are left untouched.

\subsection{Chromatic polynomials}

For a graph $G$ on the vertex set $[n]$, the chromatic polynomial $P_G(t)$ has an expansion of the form
\[
P_G(t)=\sum_{k=0}^{n}(-1)^{n-k}a_k t^k
\]
and, by a recent result of Huh~\cite{Huh}, we know the sequence $(a_0,a_1,\dots,a_n)$ is log-concave. The proof is a high-powered one using ideas from algebraic geometry, and is very far from a combinatorial proof via an explicit injection. Note that necessary and sufficient conditions (the existence of a perfect elimination order) for $(-1)^n P_G(-t)={\rm ISF}(\bbone,t)$ was given by Hallam and Sagan in~\cite{HallamS}. A refinement involving the notion of tight forest was also derived~\cite[Theorem 5.11]{HallamMS}. By Whitney's no broken circuit theorem which relies on the choice of an ordering of the edges of $G$, the coefficient $a_k$ counts all forests in $G$ with $k$ components and which do not contain a broken circuit. The latter is a set of edges obtained by removing from a circuit its minimal edge. Again on the vertex set $[n]$ with its usual ordering, let us order edges $(i,j)$ with $i<j$ using the lexicographic order, and let us call a forest admissible if it contains no broken circuit for this choice of edge ordering. It is not difficult to reformulate admissibility in a way similar to the increasing property, as follows. As before let us pick the minimal vertex as the root in each tree of the forest. For a vertex $v$, let $C(v)$ denote the set of children of $v$, and let $B(v)$ denote the branch at $v$, i.e., the union of $\{v\}$, $C(v)$, the set of grandchildren of $v$, etc. We will call a vertex $v$ a good vertex if for all $w$ in $C(v)$, $w$ is the smallest element in $B(w)$ with an edge in $G$ connecting it to $v$. A forest is admissible if and only if all the vertices are good. Note that a proof of the no broken circuit theorem with this particular edge ordering and notion of admissible forests, using the combinatorics of perturbative renormalization in quantum field theory (featuring toy versions of the Bogoliubov induction and Zimmermann's forest formula) was given in~\cite[Ch. 4]{Abdesselam}.
When $G$ is a tree, our injective proof not only applies 
but degenerates into the easier case from \S\ref{subsetsec} for subsets.
When $G$ is cycle, which is slightly less trivial, one can also find movable edges. However, if $G$ is say a triangle and an edge, as in
\[
G=\{(1,4),(2,4),(2,3),(3,4)\}
\]
with $n=4$, then movable edges may not exist, let alone be used to construct an injection.
For instance, take the admissible forests
\[
\mathscr{A}=\{(1,4),(2,4),(3,4)\}
\]
and
\[
\mathscr{B}=\{(2,3),(3,4)\}\ ,
\]
then there is no movable adge. Namely, there is no $e\in\mathscr{A}\backslash\mathscr{B}$ such that $\mathscr{A}':=\mathscr{A}\backslash\{e\}$ and
$\mathscr{B}':=\mathscr{B}\cup\{e\}$
are admissible forests. However if one changes the vertex ordering, or equivalently relabels the vertices so now
\[
G=\{(1,2),(2,3),(2,4),(3,4)\}\ ,
\]
this problem does not occur. Echoing~\cite{Pak}, we believe that a combinatorial proof of Huh's theorem is still out of reach. Yet, it may be interesting to investigate the class of graphs for which, modulo a suitable choice of vertex ordering, movable edges always exist.

\bigskip
\noindent{\bf Acknowledgements:}
{\small
For useful correspondence, the author thanks Christian Krattenthaler, Igor Pak, Bruce Sagan and Jiang Zeng.
The author also thanks Ken Ono for introducing him to a related problem on D'Arcais polynomials, which through many meanders and detours, led to the investigation in the present article.
}

\end{document}